\documentclass[11pt,reqno]{amsart}
\usepackage[T1]{fontenc}
\usepackage{lmodern,amssymb,mathtools}
\usepackage[hidelinks,pdfauthor={Liangyu Zhu},pdftitle={A Note on the Cosupport of Balmer-Favi Idempotents}]{hyperref}
\usepackage[hmargin=1.02in,vmargin=0.98in]{geometry}
\newtheorem{theorem}{Theorem}[section]
\newtheorem{proposition}[theorem]{Proposition}
\newtheorem{corollary}[theorem]{Corollary}
\newtheorem{lemma}[theorem]{Lemma}
\theoremstyle{definition}
\newtheorem{convention}[theorem]{Convention}
\newtheorem{remark}[theorem]{Remark}
\newtheorem*{mainA}{Theorem A}
\newtheorem*{mainB}{Theorem B}

\DeclareMathOperator{\Cosupp}{Cosupp}
\DeclareMathOperator{\Supp}{Supp}
\DeclareMathOperator{\gen}{gen}
\DeclareMathOperator{\Spc}{Spc}
\DeclareMathOperator{\Spec}{Spec}
\DeclareMathOperator{\Proj}{Proj}
\DeclareMathOperator{\Hom}{Hom}
\DeclareMathOperator{\Ext}{Ext}
\DeclareMathOperator{\StMod}{StMod}
\DeclareMathOperator{\fib}{fib}
\newcommand{\T}{\mathcal T}
\newcommand{\U}{\mathcal U}
\newcommand{\D}{\mathrm D}
\newcommand{\one}{\mathbf 1}
\newcommand{\p}{\mathfrak p}
\newcommand{\q}{\mathfrak q}
\newcommand{\m}{\mathfrak m}
\newcommand{\RHom}{\mathbf R\!\operatorname{Hom}}
\newcommand{\RG}{\mathbf R\Gamma}
\newcommand{\dt}{\otimes^{\mathbf L}}
\newcommand{\fg}{\mathrm{fg}}
\newcommand{\fl}{\mathrm{fl}}
\newcommand{\qc}{\mathrm{qc}}
\newcommand{\kpp}{\kappa(\p)}
\numberwithin{equation}{section}
\title[Cosupport of Balmer--Favi idempotents]{A Note on the Cosupport of Balmer--Favi Idempotents}
\author{Liangyu Zhu}
\address{School of Mathematical Sciences, Anhui University, Hefei 230601, China}
\email{a126112006@stu.ahu.edu.cn}
\date{14 November 2026}
\keywords{Tensor triangular geometry, cosupport, Balmer--Favi idempotent, local duality, stable module category, monochromatic layer}
\subjclass[2020]{Primary 18F99; Secondary 13D45, 20J06, 55P42}
\begin{document}
\begin{abstract}
Let $g_x$ be the Balmer--Favi idempotent at a point $x$ of the spectrum of a rigidly-compactly generated tensor triangulated category. We study the equality $\Cosupp(g_x)=\gen(x)$ by expressing $g_x$ as a local dual. We prove the equality for the derived category of a Noetherian scheme and for the stable module category of a finite group scheme over a field of positive characteristic. For $E(n)$-local spectra we express the maps between distinct monochromatic layers of the sphere in terms of the chromatic fracture square. The height-two splitting gives the equality through height two at every prime. For a Noetherian local ring we characterize the Gorenstein property by the form of its closed-point idempotent.
\end{abstract}
\maketitle
\section{Introduction and main results}

For a commutative Noetherian ring $R$, the support of a perfect complex $c$ consists of the primes $\p$ for which $c_{\p}\ne0$. In particular, support determines the thick tensor ideal generated by $c$. Balmer's spectrum extends this description to an essentially small rigid tensor triangulated category $\mathcal K$: a point of $\Spc(\mathcal K)$ is a prime thick tensor ideal $\mathcal P$, and $\mathcal P$ belongs to the support of $c$ precisely when $c\notin\mathcal P$ \cite{bal}. For $\mathcal K=\D^{\mathrm{perf}}(R)$ the prime corresponding to $\p$ consists of the perfect complexes which vanish over $R_{\p}$.

Let $\T$ be a rigidly-compactly generated tt-category with $\mathcal K=\T^c$ and weakly Noetherian spectrum. For each point $x$, the construction of Balmer and Favi gives a tensor idempotent $g_x$ \cite{bf,bchs}. Since $g_x\otimes g_x\simeq g_x\ne0$ and $g_y\otimes g_x=0$ for $y\ne x$, the Balmer--Favi support of $g_x$ is the singleton $\{x\}$.

Write $[a,b]$ for the internal Hom, that is, the object representing $t\mapsto\Hom_{\T}(t\otimes a,b)$. The cosupport of $b$ consists of the points $y$ for which $[g_y,b]\ne0$ \cite{bchs}. The identity of $g_x$ gives $x\in\Cosupp(g_x)$. Barthel, Castellana, Heard and Sanders prove the other inclusion in
\begin{equation}\label{e:bound}
 \{x\}\subseteq\Cosupp(g_x)\subseteq\gen(x),
 \qquad \gen(x)=\{y\mid x\in\overline{\{y\}}\}
\end{equation}
\cite[Lemma~4.26]{bchs}. They ask for which categories the right-hand inclusion is an equality \cite[Question~21.3]{bchs}. Thus the question is a nonvanishing problem for each pair of points in the specialization order:
\begin{equation}\label{e:question}
 x\in\overline{\{y\}}\quad\Longrightarrow\quad [g_y,g_x]\ne0\;?
\end{equation}

Let $R$ be a discrete valuation ring, $\pi$ a uniformizer, $K$ the fraction field, $\eta$ the generic point and $s$ the closed point. The local cohomology triangle for $(\pi)$ is
\[
 g_s\longrightarrow R\longrightarrow K\longrightarrow g_s[1].
\]
Here $g_\eta=K$ and $g_s\simeq(K/R)[-1]$. The quotient map $K\to K/R$ gives a nonzero morphism $g_\eta\to g_s[1]$. Hence $\eta\in\Cosupp(g_s)$, although $K\otimes_R^{\mathbf L}g_s=0$. More generally, \cite[Example~8.16]{bchs} establishes \eqref{e:question} for every big tt-category whose spectrum consists of a closed point and its unique strict generalization. For a three-point chain with $z\in\overline{\{y\}}$ and $y\in\overline{\{x\}}$, \eqref{e:question} also requires $[g_x,g_z]\ne0$. Nonzero morphisms between successive layers may have zero composite.

Let $(B,\m,k)$ be local, let $E=E_B(k)$, and suppose first that $B$ has a normalized dualizing complex $\omega_B^\bullet$. Local duality identifies $g_{\m}$ with $\RHom_B(\omega_B^\bullet,E)$. For a prime $\q\subseteq\m$ the tensor--Hom adjunction gives
\[
 [g_{\q},g_{\m}]
 \simeq\RHom_B(g_{\q}\otimes_B^{\mathbf L}\omega_B^\bullet,E).
\]
The complex on the right is nonzero: $\omega_B^\bullet$ has support at $\q$, and $\Hom_B(-,E)$ detects nonzero $B$-modules. Hence $\Cosupp_B(g_{\m})=\Spec B$.

Let $A$ be a graded-commutative Noetherian ring acting through the graded endomorphisms of the tensor unit. For a homogeneous prime $\p$, let $I(\p)$ be the graded injective hull of $A/\p$. Following Benson, Iyengar and Krause, let $J_{\p}$ represent the functor $z\mapsto\Hom_A(\Hom^*_{\T}(\one,z),I(\p))$, where $\Hom_A$ denotes degree-zero graded maps \cite[(5.3)]{bik}. Their Proposition~5.4 computes the cosupport of $T_c(I(\p))$ for compact $c$. Rigidity gives $T_c(I(\p))\simeq[c^\vee,J_{\p}]$. We allow an arbitrary $W$ in $[W,J_{\p}]$.

\begin{mainA} Under Convention~\ref{conv}, for each $W\in\T$ one has
\[
 \Cosupp([W,J_{\p}])=\Supp(W)\cap\gen(\p).
\]
Consequently, if $g_{\p}$ is a shift of $[W,J_{\p}]$, it has cosupport $\gen(\p)$ precisely when $\Supp(W)$ contains all generalizations of $\p$.
\end{mainA}

Hrbek, Nakamura and \v{S}\v{t}ov\'{\i}\v{c}ek compute $\RHom_R(g_{\q},g_{\p})$ for a saturated inclusion $\q\subsetneq\p$ \cite[Theorem~6.1]{hns}. For any inclusion $\q\subseteq\p$, their Proposition~6.9 describes the vanishing of its cohomology in terms of a dualizing complex over $\widehat{R_{\p}}$. We give the relation with \eqref{e:question} after Corollary~\ref{c:sch}.

For an arbitrary Noetherian ring $R$, the completion of $R_{\p}$ admits a dualizing complex. Restriction of scalars sends its closed-point idempotent to $g_{\p}$ (Lemma~\ref{l:comp} and Corollary~\ref{c:sch}). For finite group schemes, local duality identifies $\Gamma_{\p}k$ with a suspension of the dual of the modular character with respect to $J_{\p}$ \cite[Theorem~5.1]{bikp}. For $E(n)$-local spectra the chromatic fracture square gives the function spectrum in Corollary~\ref{c:chrom}.

\begin{mainB}
The equality $\Cosupp(g_x)=\gen(x)$ holds at every point of the following categories:
\begin{enumerate}
\item $\D_{\qc}(X)$ for a Noetherian scheme $X$;
\item $\StMod(G)$ for a finite group scheme $G$ over a field of positive characteristic;
\item the $E(n)$-local stable homotopy category at any prime for $n\le 2$.
\end{enumerate}
More generally, in the $E(n)$-local category the equality holds at the points of height at most two for every $n$.
\end{mainB}

For a Cohen--Macaulay local ring with canonical module $\omega_B$, the normalized dualizing complex is $\omega_B[\dim B]$. Consequently $g_{\m}$ is the Matlis dual of $\omega_B$, placed in cohomological degree $\dim B$. We characterize the primes at which this object is a shift of an indecomposable injective, and give the duality formula for complexes with finite-length cohomology (Corollaries~\ref{c:cm}, \ref{c:finite} and Theorem~\ref{t:bass}). 

\section*{Acknowledgments}
The author thanks Greg Stevenson for helpful discussions in February 2026. The author also thanks his master's advisor, Hongxing Chen, for discussions in March 2026 on the commutative algebra in this paper.

\section{Point idempotents and local duality}

We write $[a,b]$ for the internal Hom in $\T$. For a Thomason subset $V$ of $\Spc(\T^c)$, let $\mathbb{E}_V\to\one\to \mathbb{F}_V\to\Sigma \mathbb{E}_V$ denote the associated idempotent triangle. A point $x$ is weakly visible if $\{x\}=V\cap W^c$ for two Thomason subsets. In this case $g_x=\mathbb{E}_V\otimes \mathbb{F}_W$ is independent of this presentation up to isomorphism. Weak Noetherianity means that all points are weakly visible. We use the conventions
\[
 \Supp(a)=\{x\mid g_x\otimes a\ne0\},\qquad
 \Cosupp(b)=\{x\mid[g_x,b]\ne0\}.
\]
For a compact object $c$, write $c^\vee=[c,\one]$. Our convention is $\Hom^n_{\T}(a,b)=\Hom_{\T}(a,\Sigma^n b)$, and $\Hom^*_{\T}(a,b)$ is the direct sum over $n\in\mathbb Z$. For a graded module $N$ we use $N(d)^i=N^{i+d}$. A Thomason subset is a union of closed subsets with quasi-compact open complements. For $\mathbb{E}_V$ the functor $\mathbb{E}_V\otimes-$ is the acyclization associated to $V$, and $\mathbb{F}_V\otimes-$ is the complementary smashing localization.

\begin{convention}\label{conv}
Let $A$ be a graded-commutative Noetherian ring acting through the graded endomorphisms of $\one$. Suppose the spectrum of compact objects is identified with a subspace $P$ of $\Spec^{\mathrm h}A$, with its induced topology. We use homogeneous localization. Suppose that for $\q\in P$, $c\in\T^c$ and $a\in\T$ the graded $A$-module
\[
 \Hom^*_{\T}(c,g_{\q}\otimes a)
\]
is $\q$-local and $\q$-torsion. Thus homogeneous elements outside $\q$ act invertibly, and each element of the module is killed by a power of $\q$. These conditions hold for the point local cohomology functors for the canonical actions on $\D(R)$ and $\StMod(G)$; see \cite{bik,bikp}.
\end{convention}

Let $I(\p)$ denote the graded injective envelope of $A/\p$. Brown representability gives an object $J_{\p}=T_{\one}(I(\p))$ with natural isomorphisms
\begin{equation}\label{e:brown}
 \Hom_{\T}^*(z,J_{\p})\cong
 \underline{\Hom}_A^*(\Hom^*_{\T}(\one,z),I(\p)).
\end{equation}
Here the right-hand side denotes graded Hom in all degrees. Its degree-zero part defines a cohomological functor taking coproducts to products, since $\one$ is compact and $I(\p)$ is injective. This is the representing object of \cite[(5.3)]{bik}.

\begin{lemma}\label{l:test}
Let $N$ be a $\q$-local, $\q$-torsion graded $A$-module. Then
$\underline{\Hom}^*_A(N,I(\p))\ne0$ precisely when $N\ne0$ and $\q\subseteq\p$.
\end{lemma}
\begin{proof}
If $\q\subseteq\p$, every homogeneous element outside $\p$ acts invertibly on $N$. Thus $N$ is an $A_{\p}$-module. The shifts of $I(\p)$ cogenerate the category of graded $A_{\p}$-modules \cite[proof of Proposition~5.4]{bik}. It follows that $\underline{\Hom}_A^*(N,I(\p))$ detects whether $N$ is zero.

If $\q\not\subseteq\p$, choose a homogeneous $s\in\q\setminus\p$. For a homogeneous map $h:N\to I(\p)(a)$ and $v\in N$, choose $t$ with $s^tv=0$. Then $s^th(v)=0$. Since $s$ acts invertibly on $I(\p)$, we get $h(v)=0$.
\end{proof}

\begin{theorem}\label{t:pairing}
Under Convention~\ref{conv}, for every $W\in\T$ and $\p\in P$,
\begin{equation}\label{e:pairing}
 \Cosupp([W,J_{\p}])=\Supp(W)\cap\gen(\p).
\end{equation}
In particular, suppose $g_{\p}\simeq\Sigma^a[W,J_{\p}]$ for some integer $a$. Then $\Cosupp(g_{\p})=\gen(\p)$ if and only if $\gen(\p)\subseteq\Supp(W)$.
\end{theorem}
\begin{proof}
Fix $\q\in P$ and put $Z=g_{\q}\otimes W$. We must decide whether $[g_{\q},[W,J_{\p}]]\simeq[Z,J_{\p}]$ vanishes. Test this object against $c^\vee$ for $c\in\T^c$. We have
\begin{align*}
 \Hom^*_{\T}(c^\vee,[Z,J_{\p}])
 &\cong \Hom^*_{\T}(c^\vee\otimes Z,J_{\p})\\
 &\cong \underline{\Hom}^*_A
       \big(\Hom_{\T}^*(\one,c^\vee\otimes Z),I(\p)\big)\\
 &\cong \underline{\Hom}^*_A
       \big(\Hom_{\T}^*(c,Z),I(\p)\big).
\end{align*}
The second line is \eqref{e:brown}. For the last line use $c^\vee\otimes Z\simeq[c,Z]$. Convention~\ref{conv} applies to $\Hom_{\T}^*(c,Z)$ for every $W$.

If $\q\not\subseteq\p$, Lemma~\ref{l:test} makes the last group zero for every $c$. The $c^\vee$ range over all compact objects, so $[Z,J_{\p}]=0$. If $\q\subseteq\p$, the last group is nonzero precisely when $\Hom^*_{\T}(c,Z)$ is nonzero. Such a $c$ exists precisely when $Z\ne0$, by compact generation. Thus
\[
 \q\in\Cosupp([W,J_{\p}])
 \quad\Longleftrightarrow\quad \q\subseteq\p\ \text{and}\ \q\in\Supp(W).
\]
Since $\gen(\p)=\{\q\in P\mid\q\subseteq\p\}$, this proves \eqref{e:pairing}. Taking a suspension leaves cosupport unchanged.
\end{proof}

\begin{remark}
The case $W=\one$ recovers the unit case of \cite[Proposition~5.4]{bik} when the central and tensor triangular point functors agree. If $u$ is tensor-invertible, $[u,J_{\p}]\simeq u^{-1}\otimes J_{\p}$ and $\Supp(u)=P$.
\end{remark}

\begin{proposition}\label{p:transport}
Let $f^*:\T\to\U$ be a geometric functor between rigidly-compactly generated tt-categories with weakly Noetherian spectra. Write $f_*$ for its right adjoint and $\varphi$ for the induced map on spectra. Suppose $f_*$ is conservative and $\U$ satisfies codetection. If $f_*g_u\simeq g_x$ and $\varphi(\gen(u))=\gen(x)$, then
\[
 \Cosupp_{\U}(g_u)=\gen(u)\quad\Longrightarrow\quad
 \Cosupp_{\T}(g_x)=\gen(x).
\]
More generally, $\Cosupp_{\T}(f_*b)=\varphi(\Cosupp_{\U}(b))$ for each $b\in\U$. Under detection in $\U$, the analogous equality holds for support. For arbitrary $f_*$, codetection gives the cosupport inclusion $\subseteq$, and detection gives the support inclusion $\subseteq$.
\end{proposition}
\begin{proof}
The base-change formulas for $f_*b$, including the inclusions for arbitrary $f_*$, are \cite[Proposition~14.2]{bchs}. For $b=g_u$ the cosupport formula gives
\[
 \Cosupp_{\T}(g_x)=\varphi(\Cosupp_{\U}(g_u))
                =\varphi(\gen(u))=\gen(x).\qedhere
\]
\end{proof}

\begin{corollary}\label{c:model}
Under the hypotheses of Proposition~\ref{p:transport}, suppose that $\U$ satisfies Convention~\ref{conv}, $u$ corresponds to $\p$, and $g_u\simeq\Sigma^a[W,J_{\p}]$. Then
\[
 \Cosupp_{\T}(g_x)
   =\varphi\big(\Supp_{\U}(W)\cap\gen(u)\big).
\]
If $\gen(u)\subseteq\Supp_{\U}(W)$, then $\Cosupp_{\T}(g_x)=\gen(x)$. The same conclusion holds for a composite of such changes of category whenever the successive right adjoints carry the indicated point idempotents to one another.
\end{corollary}
\begin{proof}
Proposition~\ref{p:transport} and \eqref{e:pairing} give
\[
 \Cosupp_{\T}(g_x)=\varphi(\Cosupp_{\U}(g_u))
 =\varphi\big(\Supp_{\U}(W)\cap\gen(u)\big).
\]
For composable geometric functors the right adjoints compose in reverse order and the induced maps on spectra compose in the same order as the right adjoints.
\end{proof}

\begin{proposition}\label{p:unit}
Let $f^*:\T\to\U$ be as in Proposition~\ref{p:transport}, and let $x,u$ be closed points with $\varphi^{-1}(x)=\{u\}$ and both singletons Thomason. Put $\mathbb{E}=g_x$ and $\mathbb{E}'=g_u$. Suppose $f_*$ preserves coproducts. To obtain $\mathbb{E}\simeq f_*\mathbb{E}'$ it suffices that the unit $a\to f_*f^*a$ be an isomorphism on a set of localizing generators for $\mathbb{E}\otimes\T$. It also suffices that it be an isomorphism on $\mathbb{E}$ alone.
\end{proposition}
\begin{proof}
Let $\eta:1_{\T}\to f_*f^*$ be the adjunction unit. Since $f^*$ and $f_*$ are exact and preserve coproducts, $\{a\mid\eta_a\text{ is invertible}\}$ is a localizing subcategory. If it contains the given generators, it contains $\mathbb E\in\mathbb E\otimes\T$. The inverse-image formula gives $f^*\mathbb E\simeq\mathbb E_{\varphi^{-1}\{x\}}=\mathbb E'$, so $\eta_{\mathbb E}$ identifies $\mathbb E$ with $f_*\mathbb E'$.
\end{proof}

\begin{remark}\label{r:hom}
In a category satisfying the Hom formula $\Cosupp([a,b])=\Supp(a)\cap\Cosupp(b)$, any isomorphism $g_u\simeq[W,I]$ yields $\Cosupp(g_u)=\Supp(W)\cap\Cosupp(I)$. In a stratified category with weakly Noetherian spectrum both this formula and the tensor product formula hold; see \cite[Theorem~7.15 and Remark~14.23]{bchs}. In particular $\Supp(g_x)=\{x\}$ by orthogonality and nonvanishing of the point idempotents. These formulas apply to $\D(R)$ for a commutative Noetherian ring $R$ \cite[Example~19.2]{bchs}.
\end{remark}

\section{Local duality and the closed-point idempotent}

In a local category let $x$ be a closed point such that $\{x\}$ is Thomason, and write $\mathbb{E}=\mathbb{E}_{\{x\}}=g_x$. It comes with a counit $\mathbb{E}\to\one$. Write $\Gamma=\mathbb{E}\otimes-$ and call an object $z$ torsion at $x$ if $\Gamma z\to z$ is an isomorphism. 

\begin{theorem}\label{t:duals}
Let $D\in\T$, put $I=\Gamma D$, and write $z^\dagger=[z,D]$. Let $\mathcal C$ be a full subcategory closed under $(-)^\dagger$. Suppose the biduality map $X\to[X^\dagger,D]$ is invertible for each $X\in\mathcal C$. Suppose that for $Y$ in $\mathcal C$ the canonical map
\begin{equation}\label{e:comparison}
 \mathbb{E}\otimes[Y,D]\longrightarrow [Y,\mathbb{E}\otimes D]
\end{equation}
is invertible. For $X\in\mathcal C$ there is then a natural isomorphism
\begin{equation}\label{e:dual}
 \Gamma X\simeq[X^\dagger,I].
\end{equation}
If $\one\in\mathcal C$, then $\mathbb{E}\simeq[D,I]$. If also $D=\Sigma^d W$, then $\mathbb{E}\simeq\Sigma^{-d}[W,I]$ and $\Gamma W\simeq\Sigma^{-d} I$. In particular $D\simeq\Sigma^d\one$ gives $\mathbb{E}\simeq\Sigma^{-d}I$.
\end{theorem}
\begin{proof}
Let $\delta_X:X\to[X^\dagger,D]$ be the biduality map. Write $\theta_Y$ for the adjoint of
\[
 (\mathbb E\otimes[Y,D])\otimes Y
 \xrightarrow{1\otimes\mathrm{ev}}\mathbb E\otimes D.
\]
Since $X^\dagger\in\mathcal C$, the following composite consists of natural isomorphisms:
\[
 \Gamma X\xrightarrow{1\otimes\delta_X}\mathbb E\otimes[X^\dagger,D]
 \xrightarrow{\theta_{X^\dagger}}[X^\dagger,\mathbb E\otimes D]=[X^\dagger,I].
\]
At $X=\one$ this gives $\mathbb E\simeq[D,I]$. For $D=W[d]$, one has $[W[d],I]\simeq[W,I][-d]$ and $\Gamma W[d]\simeq I$.
\end{proof}

\begin{remark}
For dualizable $Y$, \eqref{e:comparison} is an isomorphism since both sides identify with $Y^\vee\otimes \mathbb{E}\otimes D$. Lemma~\ref{l:input} establishes the comparison for bounded complexes with finite cohomology over a Noetherian local ring with a normalized dualizing complex.
\end{remark}

\begin{proposition}\label{p:length}
Suppose the hypotheses of Theorem~\ref{t:duals} hold in a category with a $t$-structure. Assume $D=W[d]$, and assume there is an exact contravariant functor $D_0$ on the heart and natural isomorphisms
$H^i([z,I])\cong D_0(H^{-i}z)$ for every $z=X^\dagger$ with $X\in\mathcal C$. Then
\begin{equation}\label{e:finite}
 H^i(\Gamma X)\cong D_0\big(H^{d-i}([X,W])\big).
\end{equation}
For $X$ torsion at $x$, the left side is $H^i(X)$. Suppose further that $\mathcal C$ is closed under truncations and that its objects are bounded. If $\Gamma$ fixes the objects of a Serre subcategory $\mathcal A$ of the heart, then \eqref{e:finite} computes $H^i(X)$ for every $X\in\mathcal C$ whose cohomology belongs to $\mathcal A$.

For an object $N$ of $\mathcal C$ in the heart, assume $D_0$ detects zero on the objects $H^j([N,W])$ for all $j$. Then
\[
 H^j([N,W])=0\ (j\ne d)
 \quad\Longleftrightarrow\quad
 H^i(\Gamma N)=0\ (i\ne0).
\]
\end{proposition}
\begin{proof}
For $X\in\mathcal C$, $X^\dagger=[X,W][d]$. Hence \eqref{e:dual} gives
\[
 H^i(\Gamma X)\cong H^i([X^\dagger,I])
 \cong D_0(H^{-i}X^\dagger)=D_0(H^{d-i}([X,W])).
\]
Suppose $X\in\T^{\ge a}\cap\T^{\le b}$ and $H^j(X)\in\mathcal A$ for all $j$. The counit is invertible at $\tau^{\le a-1}X=0$ and at each $H^j(X)[-j]$. Apply the counit to the truncation triangle
\[
 \tau^{\le j-1}X\longrightarrow\tau^{\le j}X
 \longrightarrow H^j(X)[-j]\longrightarrow(\tau^{\le j-1}X)[1].
\]
Induction on $j$ and the two-out-of-three property give $\Gamma\tau^{\le b}X\simeq\tau^{\le b}X=X$. For $N$ in the heart, \eqref{e:finite} and detection of zero by $D_0$ give $H^i(\Gamma N)=0\iff H^{d-i}([N,W])=0$. Substitute $j=d-i$.
\end{proof}

\begin{proposition}\label{p:residue}
Suppose $\one\in\mathcal C$ in Theorem~\ref{t:duals}. Assume the canonical maps $z\to[[z,I],I]$ are isomorphisms for $z=D,\one$. Then
\[
 \mathbb{E}\simeq I[-d]\quad\Longleftrightarrow\quad D\simeq\one[d].
\]
In either case there is an isomorphism $[t,\one]\simeq[t,I][-d]$ for every torsion object $t$.
\end{proposition}
\begin{proof}
Theorem~\ref{t:duals} gives $\mathbb E\simeq[D,I]$. If $\mathbb E\simeq I[-d]$, then biduality at $D$ and $\one$ gives
\[
 D\simeq[[D,I],I]\simeq[I[-d],I]\simeq[I,I][d]\simeq\one[d].
\]
Conversely $D\simeq\one[d]$ gives $\mathbb E\simeq[\one[d],I]\simeq I[-d]$. For $\Gamma t\simeq t$, one has $\Gamma(a\otimes t)\simeq a\otimes t$ for every $a$. The colocalization counit therefore induces the middle isomorphism in
\[
 \Hom_{\T}(a,[t,\mathbb E])\cong\Hom_{\T}(a\otimes t,\Gamma\one)
 \xrightarrow{\sim}\Hom_{\T}(a\otimes t,\one)\cong\Hom_{\T}(a,[t,\one]).
\]
Yoneda gives $[t,\mathbb E]\simeq[t,\one]$ for every such $t$. Under the stated equivalent conditions this is $[t,I][-d]\simeq[t,\one]$.
\end{proof}

\section{Fracture triangles and generalizations}

\begin{proposition}\label{p:fract}
Suppose there is a triangle $g_x\to b\to h\to\Sigma g_x$ in $\T$. If $[g_y,b]=0$, then
\[
 [g_y,g_x]\simeq\Sigma^{-1}[g_y,h].
\]
Under this vanishing hypothesis, if a retract $r$ of $h$ admits a map $g_y\to\Sigma^a r$ which becomes nonzero under an exact functor, then $y\in\Cosupp(g_x)$. The same conclusion holds if some suspension of $g_y$ is a retract of $h$.
\end{proposition}
\begin{proof}
Applying $[g_y,-]$ gives the triangle
\[
 [g_y,g_x]\longrightarrow0\longrightarrow[g_y,h]
 \xrightarrow{\partial}\Sigma[g_y,g_x],
\]
so $\partial$ is invertible. Write $j:r\to h$, $p:h\to r$ with $pj=1_r$. For $0\ne v:g_y\to\Sigma^a r$ one has $(\Sigma^ap)(\Sigma^aj)v=v\ne0$. Thus $\Hom_{\T}(\one,\Sigma^a[g_y,h])\ne0$ by adjunction, and $[g_y,g_x]\ne0$. A morphism with nonzero image under an exact functor is nonzero. For $r\simeq\Sigma^b g_y$, take $a=-b$ and $v$ an isomorphism.
\end{proof}

\begin{corollary}\label{c:layer}
Let $L_{-1}=0,L_0,\ldots,L_n$ be nested smashing localizations and put $g_i=\fib(L_i\one\to L_{i-1}\one)$. Assume that these are the point idempotents of a chain whose generalizations of $j$ are $0,\ldots,j$. Let $Q_j$ be a localization and suppose there is a fracture triangle
\[
 g_j\longrightarrow Q_j\one\longrightarrow L_{j-1}Q_j\one\longrightarrow\Sigma g_j
\]
and suppose $Q_j(a\otimes g_i)=0$ for every $a\in\T$ and $i<j$. Then $\Cosupp(g_j)=\{0,\ldots,j\}$ if and only if $[g_i,L_{j-1}Q_j\one]\ne0$ for all $i<j$. It suffices, for each $i<j$, to exhibit a retract as in Proposition~\ref{p:fract}.
\end{corollary}
\begin{proof}
For $i<j$ and every $a\in\T$, localness of $Q_j\one$ gives
\[
 \Hom_{\T}(a,[g_i,Q_j\one])
 \cong\Hom_{\T}(Q_j(a\otimes g_i),Q_j\one)=0.
\]
Hence $[g_i,Q_j\one]=0$. Proposition~\ref{p:fract} identifies $[g_i,g_j]$ with $\Sigma^{-1}[g_i,L_{j-1}Q_j\one]$. By \eqref{e:bound}, $j\in\Cosupp(g_j)\subseteq\{0,\ldots,j\}$, proving the equivalence.
\end{proof}

\section{Examples}\label{s:examples}
\subsection{Commutative algebra and algebraic geometry}
Let $R$ be a commutative Noetherian ring. The identification $\Spec R\cong\Spc(\D^{\mathrm{perf}}(R))$ sends $\p$ to the prime of perfect complexes vanishing over $R_{\p}$. Under this identification,
\[
 g_{\p}=\RG_{\p R_{\p}}R_{\p},\qquad
 J_{\p}=E_R(R/\p),\qquad
 \Supp(X^\bullet)=\{\q\mid\kappa(\q)\dt_R X^\bullet\ne0\}.
\]
One first localizes at $\p$, then takes local cohomology at the closed point of $\Spec R_{\p}$. On $R_{\p}$ these operations give $\RG_{\p R_{\p}}R_{\p}$, computed by the tensor product of the two-term complexes $(R_{\p}\to R_{\p}[1/s])$ for generators $s$ of $\p R_{\p}$ in degrees zero and one. The second follows from exactness of Hom into the injective hull, and the third because the residue field and the point idempotent generate the same localizing tensor ideal \cite{bik,bchs}. To check Convention~\ref{conv}, let $c$ be perfect. By rigidity the groups $\Hom^*_R(c,g_{\q}\dt_R a)$ are the cohomology of $g_{\q}\dt_R(c^\vee\dt_R a)$. This complex is obtained by taking local cohomology at $\q R_{\q}$ after localization at $\q$, so its cohomology is $\q$-local and $\q$-torsion. For a bounded-above complex with finitely generated cohomology the last support agrees with $\{\q\mid X^\bullet_{\q}\ne0\}$: if $m$ is the largest index with $H^m(X^\bullet_{\q})\ne0$, the truncation triangle and right $t$-exactness of derived tensor with $\kappa(\q)$ give
\[
 H^m(\kappa(\q)\dt_R X^\bullet)=\kappa(\q)\otimes_{R_{\q}}H^m(X^\bullet_{\q})\ne0
\]
by Nakayama's lemma. Conversely, $X^\bullet_{\q}=0$ implies $\kappa(\q)\dt_R X^\bullet=0$. In particular this applies to pseudo-coherent complexes. For arbitrary complexes $X^\bullet,Y^\bullet$ one has $\Supp(X^\bullet\dt_R Y^\bullet)=\Supp(X^\bullet)\cap\Supp(Y^\bullet)$ and
$\Cosupp(\RHom_R(X^\bullet,Y^\bullet))=\Supp(X^\bullet)\cap\Cosupp(Y^\bullet)$ by Remark~\ref{r:hom}. Taking $W=R$ in Theorem~A gives $\Cosupp(E_R(R/\p))=\gen(\p)$.

\begin{lemma}\label{l:input}
Let $(B,\m,k)$ be Noetherian local with normalized dualizing complex $\omega_B^\bullet$ and $I=E_B(k)$. Theorem~\ref{t:duals} applies to $\mathcal C=\D^b_{\fg}(B)$, with $\Gamma=\RG_{\m}$. Also $\Supp(\omega_B^\bullet)=\Spec B$.
\end{lemma}
\begin{proof}
The duality on $\D^b_{\fg}(B)$ and the identification $\RG_{\m}\omega_B^\bullet\simeq I[0]$ are \cite[Tags 0A7C and 0A82]{stacks}. To check \eqref{e:comparison}, choose $s_1,\ldots,s_r$ generating $\m$. The complex
\[
 K^\bullet=\bigotimes_{\nu=1}^r (B\longrightarrow B[1/s_\nu])
\]
represents $g_{\m}$ and has Tor-amplitude in $[0,r]$. An object $Y^\bullet\in\D^b_{\fg}(B)$ is pseudo-coherent and $\omega_B^\bullet$ is bounded. Thus \cite[Lemma~15.100.3(3)]{stacks}, applied with $K^\bullet,Y^\bullet,\omega_B^\bullet$, gives
\[
 K^\bullet\dt_B\RHom_B(Y^\bullet,\omega_B^\bullet)
 \xrightarrow{\ \sim\ }
 \RHom_B(Y^\bullet,K^\bullet\dt_B\omega_B^\bullet).
\]
This is the map in \eqref{e:comparison}.

For $\q\in\Spec B$ put $D^\bullet=(\omega_B^\bullet)_{\q}$. Its homothety gives $B_{\q}\simeq\RHom_{B_{\q}}(D^\bullet,D^\bullet)$. Thus $D^\bullet\ne0$, and the residue-field test for bounded complexes with finite cohomology gives $\q\in\Supp(\omega_B^\bullet)$.
\end{proof}

\begin{corollary}\label{c:loc}
With the assumptions of Lemma~\ref{l:input}, for $X^\bullet\in\D^b_{\fg}(B)$ one has naturally
\[
 \RG_{\m}X^\bullet\simeq\RHom_B(\RHom_B(X^\bullet,\omega_B^\bullet),I).
\]
In particular $g_{\m}\simeq\RHom_B(\omega_B^\bullet,I)$ and $\Cosupp_B(g_{\m})=\Spec B$.
\end{corollary}
\begin{proof}
Apply \eqref{e:dual} with $D=\omega_B^\bullet$ and $I=E_B(k)$. Taking $X^\bullet=B[0]$ gives $g_{\m}\simeq[\omega_B^\bullet,I]$. Since $I=J_{\m}$, Theorem~\ref{t:pairing} and Lemma~\ref{l:input} give
$\Cosupp_B(g_{\m})=\Supp(\omega_B^\bullet)\cap\gen(\m)=\Spec B$.
\end{proof}

\begin{lemma}\label{l:comp}
For any Noetherian local ring $(B,\m)$, restriction of scalars from $C=\widehat B$ sends $g_{\m C}$ to $g_{\m}$.
\end{lemma}
\begin{proof}
Let $\beta:B\to C$ be the completion map and put $\mathfrak n=\m C$. The ring $C/\m C=k$ has a single prime. Pullback of the Thomason idempotent for $\{\m\}$ therefore gives $\beta^*g_{\m}\simeq g_{\mathfrak n}$ \cite[Corollary~7.5]{bf}. It remains to identify the unit $g_{\m}\to\beta_*\beta^*g_{\m}$.

For fixed $i$ put $N=H^i_{\m}(B)$. Flatness of $C$ identifies the map on $H^i$ with $N\to C\otimes_B N$, $v\mapsto1\otimes v$. Write $N_t=(0:_N\m^t)$. Since $N$ is $\m$-power torsion, $N=\bigcup_{t\ge1}N_t$. On $N_t$ the unit is the scalar-extension map
\[
 N_t\longrightarrow (C/\m^t C)\otimes_{B/\m^t}N_t.
\]
It is invertible because $B/\m^t\xrightarrow{\sim}C/\m^tC$. Passing to the filtered colimit identifies $N$ with $C\otimes_B N$. Hence the unit is an isomorphism on every cohomology module, as required.
\end{proof}

\begin{corollary}\label{c:sch}
For any Noetherian ring $R$ one has $\Cosupp_R(g_{\p})=\gen(\p)$ for every $\p\in\Spec R$. For any Noetherian scheme $X$ one has $\Cosupp_{\D_{\qc}(X)}(g_x)=\gen(x)$ at every $x\in X$.
\end{corollary}
\begin{proof}
Fix $\p\in\Spec R$. Set $B=R_{\p}$, $\m=\p B$, $C=\widehat B$ and $\mathfrak n=\m C$. Consider
\[
 \Spec C\xrightarrow{\beta}\Spec B\xrightarrow{\alpha}\Spec R.
\]
The complete local ring $C$ admits a normalized dualizing complex \cite[Section~47.22]{stacks}. Its closed-point idempotent has cosupport $\Spec C$ by Corollary~\ref{c:loc}. Both restriction functors are conservative, and codetection holds in $\D(C)$ and $\D(B)$ \cite[Example~19.2]{bchs}. By Lemma~\ref{l:comp}, $\beta_*g_{\mathfrak n}\simeq g_{\m}$. Since $\beta$ is surjective by faithful flatness, Proposition~\ref{p:transport} yields
\[
 \Cosupp_B(g_{\m})=\beta(\Cosupp_C(g_{\mathfrak n}))=\Spec B .
\]
Now $\alpha_*g_{\m}$ is $\RG_{\p B}B$ regarded as an $R$-complex, hence is $g_{\p}$. A second application gives $\Cosupp_R(g_{\p})=\alpha(\Spec B)=\{\q\mid\q\subseteq\p\}$.

For $x\in X$, take an affine open neighborhood $j:U\hookrightarrow X$. Put $V=X\setminus U$. The localization $j^*$ has fully faithful right adjoint $\mathbf Rj_*$, with $\mathbf Rj_*j^*\simeq\mathbb F_V\otimes-$. Pullback gives $j^*g_x^X\simeq g_x^U$. Moreover $\mathbb E_V\otimes g_x^X=0$, since $x\notin V$. Tensoring the idempotent triangle for $V$ with $g_x^X$ now identifies the unit
\[
 g_x^X\xrightarrow{\ \sim\ }\mathbb F_V\otimes g_x^X
             \simeq\mathbf Rj_*g_x^U .
\]
Thus Proposition~\ref{p:transport} and the affine case give $\Cosupp_X(g_x^X)=j(\gen_U(x))$. Every generalization of $x$ belongs to $U$, so the right-hand side is $\gen_X(x)$.
\end{proof}

\begin{remark}
The affine nonvanishing also follows from \cite[Proposition~6.9, equation~(7.5)]{hns}. Indeed, in its local notation let $B=R_{\p}$, $C=\widehat B$ and let $D_C^\bullet$ be the dualizing complex used there. For $\q\subseteq\p$, choose $Q\in\Spec C$ over $\q B$. The $B$-point idempotent at $\q B$ becomes the unit after tensoring with $\kappa(Q)$, since $\kappa(Q)$ is a field extension of $\kappa(\q B)$. Consequently
\[
 \kappa(Q)\dt_C(g_{\q B}\dt_B D_C^\bullet)
 \simeq\kappa(Q)\dt_C D_C^\bullet\ne0 .
\]
The last inequality follows from full support of $D_C^\bullet$. Formula~(7.5) in \cite{hns} identifies $[g_{\q B},g_{\m}]$ with the Matlis dual over $C$ of a shift of $g_{\q B}\dt_B D_C^\bullet$. The injective hull of the residue field cogenerates $\mathrm{Mod}\,C$, so this dual is nonzero.

The localization step in Corollary~\ref{c:sch} gives $\Cosupp_R(g_{\p})$ as the image of the cosupport of the closed-point idempotent in $\D(R_{\p})$. Thus $\Cosupp_R(g_{\p})=\gen(\p)$ is equivalent to that local idempotent having full cosupport. If the ring or scheme admits a pseudo-dualizing complex, its stalks are dualizing and Corollary~\ref{c:loc} applies after normalization. This includes the case of a global dualizing complex.
\end{remark}

\begin{corollary}\label{c:cm}
Suppose $B=R_{\p}$ is Cohen--Macaulay of dimension $d$ and admits a dualizing complex. Let $\omega_B$ be its canonical module and $I=E_B(\kpp)$. Then, as $B$-complexes,
\[
 g_{\p}\simeq H^d_{\p B}(B)[-d]\simeq\Hom_B(\omega_B,I)[-d],
 \qquad \RG_{\p B}\omega_B\simeq I[-d].
\]
When $B$ is Gorenstein, $g_{\p}\simeq E_R(R/\p)[-d]$ after restriction to $R$.
\end{corollary}
\begin{proof}
By \cite[Tag 0AWS]{stacks}, $\omega_B^\bullet\simeq\omega_B[d]$. Theorem~\ref{t:duals} and injectivity of $I$ give
\[
 g_{\p}\simeq\RHom_B(\omega_B[d],I)
 \simeq\Hom_B(\omega_B,I)[-d],\qquad
 (\RG_{\p B}\omega_B)[d]\simeq\RG_{\p B}\omega_B^\bullet\simeq I .
\]
Consequently $H^i(g_{\p})=0$ for $i\ne d$ and $H^d(g_{\p})=H^d_{\p B}(B)$. Truncation gives the other expression for $g_{\p}$. If $B$ is Gorenstein then $\omega_B\simeq B$ and $E_B(\kpp)$ restricts to $E_R(R/\p)$ \cite[Tag 08Y8]{stacks}.
\end{proof}

\begin{corollary}\label{c:finite}
Let $(B,\m,k)$ be complete Cohen--Macaulay local of dimension $d$, let $X^\bullet\in\D^b_{\fg}(B)$, and write $D_0(-)=\Hom_B(-,E_B(k))$. The following conditions are equivalent:
\begin{enumerate}
\item $\Supp(X^\bullet)\subseteq\Cosupp(B)$;
\item $\Supp(X^\bullet)\subseteq\{\m\}$;
\item every $H^n(X^\bullet)$ has finite length;
\item $\Ext^j_B(H^n(X^\bullet),\omega_B)=0$ for every $n$ and $j\ne d$.
\end{enumerate}
When they hold there are natural isomorphisms
\[
 H^i(X^\bullet)\cong D_0\big(H^{d-i}(\RHom_B(X^\bullet,\omega_B))\big).
\]
For a Gorenstein ring the same assertions hold with $B$ in place of $\omega_B$.
\end{corollary}
\begin{proof}
Completeness gives $\Cosupp(B)=\{\m\}$ by \cite[Proposition~4.19 and Example~4.20]{bik}. This gives (1)$\iff$(2). Since $X^\bullet\in\D^b_{\fg}(B)$, $\Supp(X^\bullet)=\bigcup_n\Supp(H^n(X^\bullet))$. For a nonzero finite $B$-module $N$,
\[
 \Supp(N)\subseteq\{\m\}\iff\sqrt{\operatorname{ann}_B(N)}=\m
 \iff \m^tN=0\text{ for some }t\ge1
 \iff \operatorname{length}_B(N)<\infty .
\]
This gives (2)$\iff$(3), including $N=0$.

Let $N=H^n(X^\bullet)$. In Proposition~\ref{p:length} take $W=\omega_B$ and $I=E_B(k)$. Since $I$ is injective, $H^i(\RHom_B(Z^\bullet,I))=\Hom_B(H^{-i}(Z^\bullet),I)$. Thus
\[
 H^i_{\m}(N)\cong D_0(\Ext_B^{d-i}(N,\omega_B)).
\]
The Matlis dual $D_0$ is faithful and exact. If $N$ has finite length, its localization at each generator of $\m$ vanishes. The stable Koszul complex gives $\RG_{\m}N\simeq N[0]$, proving (4). Conversely (4) implies $H^i_{\m}(N)=0$ for $i\ne0$. Grothendieck's nonvanishing theorem \cite[Theorem~3.5.7]{bh} gives $H^{\dim N}_{\m}(N)\ne0$ for $N\ne0$, hence $\dim N=0$ and $N$ has finite length. This proves (4)$\Rightarrow$(3). Under (3), the truncation triangles of $X^\bullet$ give $\RG_{\m}X^\bullet\simeq X^\bullet$ by Proposition~\ref{p:length}, and \eqref{e:finite} gives the displayed isomorphism. For $B$ Gorenstein take $\omega_B=B$.
\end{proof}

\begin{theorem}\label{t:bass}
Let $R$ be Noetherian, $B=R_{\p}$, $d=\dim B$, $C=\widehat B$, $k=\kpp$, and $D_0^C=\Hom_C(-,E_C(k))$. The following are equivalent:
\begin{enumerate}
\item $B$ is Gorenstein;
\item $g_{\p}\simeq E_R(R/\p)[-d]$ in $\D(R)$;
\item $g_{\p}\simeq E_R(R/\q)[h]$ for some $\q\in\Spec R$ and $h\in\mathbb Z$;
\item for each $X^\bullet\in\D^b_{\fl}(C)$ there are isomorphisms, natural in $X^\bullet$,
\[
 H^i(X^\bullet)\cong D_0^C\big(H^{d-i}(\RHom_C(X^\bullet,C))\big)\qquad(i\in\mathbb Z).
\]
\end{enumerate}
Here $\D^b_{\fl}(C)$ denotes bounded complexes with finite-length cohomology. In particular, the non-Gorenstein locus consists of the points whose idempotent is not a shift of an indecomposable injective.
\end{theorem}
\begin{proof}
(1)$\Rightarrow$(2)$\Rightarrow$(3). Apply Corollary~\ref{c:cm} with $\omega_B=B$, then take $\q=\p$ and $h=-d$.

(3)$\Rightarrow$(1). Corollary~\ref{c:sch} and \eqref{e:pairing} with $W=R$ give $\gen(\p)=\gen(\q)$, hence $\p=\q$. Restriction $\D(B)\to\D(R)$ is fully faithful since $B\dt_R B\simeq B$. Thus (3) lifts to $\RG_{\p B}B\simeq E_B(k)[h]$ in $\D(B)$. Since $k$ is $\m$-torsion, the counit induces an isomorphism $\RHom_B(k,\RG_{\m} B)\to\RHom_B(k,B)$, by the counit argument in Proposition~\ref{p:residue}. Hence
\[
 \RHom_B(k,B)\simeq\RHom_B(k,\RG_{\p B}B)
 \simeq\Hom_B(k,E_B(k))[h]\simeq k[h].
\]
Here $\Hom_B(k,E_B(k))=\operatorname{soc}E_B(k)\cong k$. The Bass criterion \cite[Lemma~47.21.5]{stacks} now gives (1), since $\Ext_B^i(k,B)=0$ for $i\gg0$. By (1)$\Rightarrow$(2) the sole nonzero cohomology module of $g_{\p}$ occurs in degree $d$. In (3) it occurs in degree $-h$, so $h=-d$.

(1)$\Rightarrow$(4). Completion preserves the Gorenstein property and dimension, so $\omega_C=C$ and $\dim C=d$. Apply Corollary~\ref{c:finite} over $C$.

(4)$\Rightarrow$(1). Taking $X^\bullet=k[0]$ gives
\[
 D_0^C(\Ext^j_C(k,C))\cong
 \begin{cases} k& j=d,\\0& j\ne d.\end{cases}
\]
The Matlis dual detects zero, so $\Ext_C^j(k,C)=0$ for $j>d$. The Bass criterion gives that $C$ is Gorenstein. Apply \cite[Lemma~47.21.9]{stacks} to the flat local map $B\to C$, whose closed fiber is $k$, to obtain (1). Finally every indecomposable injective $R$-module is $E_R(R/\q)$ for some $\q$.
\end{proof}

\subsection{Stable homotopy theory}\label{s:chrom}
Fix a prime $p$ and work in $\mathrm{Sp}_{E(n)}$. Write $S$ for the $p$-local sphere, $S_p^\wedge$ for its $p$-completion, and $L_j=L_{E(j)}$ with $L_{-1}=0$. The spectrum of compact objects is the chain of heights $0,\ldots,n$ with generalizations of $j$ given by $0,\ldots,j$. Its point idempotents are $M_jS=\fib(L_jS\to L_{j-1}S)$ \cite[Remark~20.22]{bchs}. The function spectrum into a local target is computed in spectra.

\begin{corollary}\label{c:chrom}
For $0\le i<j\le n$ there is a natural equivalence
\[
 F(M_iS,M_jS)\simeq\Omega F(M_iS,L_{j-1}L_{K(j)}S).
\]
At every prime, $\Cosupp(M_jS)=\{0,\ldots,j\}$ for $j\le\min(2,n)$.
\end{corollary}
\begin{proof}
The chromatic fracture square \cite{hs} identifies the vertical fibers in
\[
 \begin{matrix}
 L_jS&\longrightarrow&L_{K(j)}S\\
 \big\downarrow&&\big\downarrow\\
 L_{j-1}S&\longrightarrow&L_{j-1}L_{K(j)}S .
 \end{matrix}
\]
Put $T_j=L_{j-1}L_{K(j)}S$. Thus $M_jS\to L_{K(j)}S\to T_j$ is a fiber sequence. For $i<j$, both $L_iS$ and $L_{i-1}S$ are $K(j)$-acyclic by the Bousfield classes of the smashing localizations $L_i$ and $L_{i-1}$. So is their fiber $M_iS$. Since $L_{K(j)}S$ is $K(j)$-local, $F(M_iS,L_{K(j)}S)=0$. Proposition~\ref{p:fract} gives the asserted equivalence.

For $j=0$ use \eqref{e:bound}. For $j=1$ the two-point calculation \cite[Example~8.16]{bchs} gives $F(M_0S,M_1S)\ne0$. This function spectrum is the same for every ambient height $n\ge1$, since the target is $E(1)$-local.

For $j=2$ we must find maps from suspensions of $M_0S$ and $M_1S$ into $T_2$. The height-two splittings exhibit the following retracts of $T_2$ at every prime:
\[
 r_0=\Sigma^{-3}L_0S_p^\wedge,\qquad r_1=L_1S_p^\wedge.
\]
At $p=2$ this is \cite[Theorem~2]{bgh}, at $p=3$ it is \cite[Theorem~5.11]{ghm}, and at $p\ge5$ it follows from the odd-prime splitting \cite{sy}, in the form recalled in \cite[\S1, after Conjecture~5]{bgh}.

Since $M_0S=H\mathbb Q$ and $L_0S_p^\wedge$ is rational, maps $M_0S\to L_0S_p^\wedge$ correspond to elements of $\pi_0 L_0S_p^\wedge=\mathbb Q_p$. The element $1$ therefore gives a nonzero map $\Sigma^{-3}M_0S\to r_0$.

For the height-one summand use the composite
\[
 M_1S\longrightarrow L_1S\longrightarrow L_1S_p^\wedge .
\]
The cofiber of the first arrow is $L_0S$, hence is $K(1)$-acyclic. For the second arrow let $C$ be the cofiber of $S\to S_p^\wedge$. The completion map is an equivalence modulo $p$, so $C\wedge S/p=0$ and multiplication by $p$ on $C$ is invertible. On $K(1)_*C$ the same multiplication is zero. Hence $K(1)_*C=0$. Since $L_1$ preserves $K(1)$-homology, both arrows induce isomorphisms on $K(1)_*$, with source $K(1)_*M_1S\cong K(1)_*S\ne0$. Their composite is nonzero.

The inclusions $r_i\to T_2$ admit retractions. The two maps therefore remain nonzero after composition with these inclusions, giving $F(M_iS,T_2)\ne0$ for $i=0,1$. The fiber equivalence gives $F(M_iS,M_2S)\ne0$ for both $i$. Together with \eqref{e:bound} this proves the claim.
\end{proof}

\begin{remark}
For $2<j\le n$, Corollary~\ref{c:layer} gives
\[
 \Cosupp(M_jS)=\{0,\ldots,j\}
 \iff F(M_iS,L_{j-1}L_{K(j)}S)\ne0\quad\text{for every }0\le i<j .
\]
\end{remark}

\subsection{Group representations}
Let $G$ be a finite group scheme over a field $k$ of positive characteristic. Its stable module category is rigidly-compactly generated, with compact objects represented by finite dimensional modules. Put $A=H^*(G,k)$. The spectrum is $\Proj A$, the point functor is $\Gamma_{\p}=\Gamma_{\p}k\otimes_k-$, and the cohomology of its values tested by compact objects is $\p$-local and $\p$-torsion \cite{strat,bikp}. These identifications verify Convention~\ref{conv}.

\begin{corollary}\label{c:groups}
For $\p\in\Proj A$ one has
\[
 \Cosupp_{\StMod(G)}(\Gamma_{\p} k)
 =\{\q\in\Proj A\mid \q\subseteq\p\}.
\]
In particular this computes the cosupport of the Balmer--Favi idempotents for a finite group algebra in modular characteristic.
\end{corollary}
\begin{proof}
Write $\delta_G$ for the modular character and $d=\dim(A/\p)$. Local duality \cite[Theorem~5.1]{bikp} gives $\Gamma_{\p}(\delta_G)\simeq\Sigma^dT_{\p}(k)$. To identify $T_{\p}(k)$, note that for any $G$-module $M$ the comparison
\[
 H^*(G,M)_{\p}\xrightarrow{\sim}
       \Hom^*_{\StMod(G)}(k,M)_{\p}
\]
is \cite[(2.7)]{bikp}. Since $I(\p)$ is $\p$-local, mapping either graded module into $I(\p)$ gives the same functor of $M$ before or after localization. Thus $T_{\p}(k)$ represents \eqref{e:brown} for the unit, and $T_{\p}(k)\simeq J_{\p}$.

Now $\Gamma_{\p}(\delta_G)=g_{\p}\otimes_k\delta_G$. The one-dimensional module $\delta_G$ is tensor-invertible, so
\[
 g_{\p}\simeq\Sigma^d(\delta_G^{-1}\otimes_kJ_{\p})
             \simeq\Sigma^d[\delta_G,J_{\p}].
\]
For every $\q$, tensoring $g_{\q}\otimes_k\delta_G$ with $\delta_G^{-1}$ recovers $g_{\q}\ne0$. Hence $\Supp(\delta_G)=\Proj A$, and Theorem~\ref{t:pairing} gives $\Cosupp(g_{\p})=\gen(\p)$. For a finite group the modular character is trivial \cite[\S5]{bikp}.
\end{proof}

\end{document}